\documentclass[12pt]{amsart}

\usepackage{amsmath,amssymb,latexsym,amsthm}

\newcommand{\sumlim}{\sum\limits}

\newtheorem{thm}{Theorem}[section]

\newtheorem{lem}[thm]{Lemma}

\newtheorem{cor}[thm]{Corollary}

\newtheorem{obs}[thm]{Observation}

\newtheorem{exa}[thm]{Example}

\theoremstyle{definition}
\newtheorem{defn}[thm]{Definition}

\newcommand{\eqdef}{:=}

\newtheorem{prop}[thm]{Proposition}
\newtheorem{conj}[thm]{Conjecture}
\newtheorem{clm}[thm]{Claim}
\newcommand{\een}{\end{enumerate}}
\newcommand{\blem}{\begin{lem}}
\newcommand{\elem}{\end{lem}}
\newcommand{\bcl}{\begin{clm}}
\newcommand{\ecl}{\end{clm}}
\newcommand{\bthm}{\begin{thm}}
\newcommand{\ethm}{\end{thm}}
\newcommand{\bpr}{\begin{prop}}
\newcommand{\epr}{\end{prop}}
\newcommand{\bco}{\begin{cor}}
\newcommand{\eco}{\end{cor}}
\newcommand{\bcon}{\begin{conj}}
\newcommand{\econ}{\end{conj}}
\newcommand{\bde}{\begin{defn}}
\newcommand{\ede}{\end{defn}}
\newcommand{\bex}{\begin{exa}}
\newcommand{\eexa}{\end{exa}}
\newcommand{\bobs}{\begin{obs}}
\newcommand{\eobs}{\end{obs}}
\newcommand{\bexe}{\begin{exe}}
\newcommand{\eexe}{\end{exe}}

\newcommand{\grn}{G_{r,n}}

\begin{document}

\title[Excedance number of colored permutation groups]{On the excedance number of colored permutation groups}

\author{Eli Bagno}
\address{Einstein institute of Mathematics, The Hebrew University,
Givat Ram, 91904 Jerusalem, Israel}
\email{bagnoe@math.huji.ac.il}

\author{David Garber}
\address{Einstein institute of Mathematics, The Hebrew University,
Givat Ram, 91904 Jerusalem, Israel, and, \\
Department of Sciences, Holon Academic Institute of Technology,
PO Box 305, 58102 Holon, Israel}
\email{garber@math.huji.ac.il,garber@hait.ac.il}

\date{\today}

\maketitle
\begin{abstract}
We generalize the results of Ksavrelof and Zeng about the
multidistribution of the excedance number of $S_n$ with some
natural parameters to the colored permutation group and to the
Coxeter group of type $D$. We define two different orders on these
groups which induce two different excedance numbers. Surprisingly,
in the case of the colored permutation group, we get the same
generalized formulas for both orders.
\end{abstract}

\bibliographystyle{is-alpha}

\section{Introduction}

Let $r$ and $n$ be two positive integers. The {\it colored permutation
group} $G_{r,n}$ consists of all permutations of the set
$$\Sigma=\{1,\dots,n,\bar{1},\dots,\bar{n},\dots,1^{[r-1]},\dots,n^{[r-1]}\}$$
satisfying $\pi(\bar{i})=\overline{\pi(i)}$.

The symmetric group $S_n$ is a special case of $\grn$ for $r=1$.
In $S_n$ one can define the following well-known parameters: Given
$\sigma \in S_n$, $i \in [n]$ is {\it an excedance of $\sigma$} if and
only if $\sigma(i)>i$. The number of excedances is denoted by
${\rm exc}(\sigma)$. Two other natural parameters on $S_n$ are the
number of fixed points and the number of cycles of $\sigma$,
denoted  by ${\rm fix}(\sigma)$ and ${\rm cyc}(\sigma)$ respectively.

Consider the following generating function over $S_n$:

$$P_n(q,t,s)=\sumlim_{\sigma \in
S_n}{q^{{\rm exc}(\sigma)}t^{{\rm fix}(\sigma)}s^{{\rm cyc}(\sigma)}}.$$

$P_n(q,1,1)$ is the classical Eulerian polynomial, while
$P_n(q,0,1)$ is the counter part for the derangements, i.e. the
permutations without fixed points, see \cite{Sta}.

In the case $s=-1$, the two polynomials $P_n(q,1,-1)$ and $P_n(q,0,-1)$
have simple closed formulas:

\begin{equation} \label{q=1}
P_n(q,1,-1)=-(q-1)^{n-1}
\end{equation}

\begin{equation}\label{q=-1}
P_n(q,0,-1)=-q[n-1]_q.
\end{equation}

Recently, Ksavrelof and Zeng \cite{KZ} proved some new
recursive formulas which induce Equations (\ref{q=1}) and (\ref{q=-1}).\\
A natural problem is to generalize the results of \cite{KZ}
to the colored permutation groups. The main challenge here is to
choose a suitable order on the alphabet $\Sigma$ of the group $\grn$ and define the parameters
properly.

In this paper we cope with this challenge. We define two different
orders on $\Sigma$, one of them 'forgets' the colors, while the
other is much more natural, since it takes into account the color
structure of $\grn$. The parameter ${\rm exc}$ will be defined
according to both orders in two different ways. The interesting
point is that for the group $\grn$ we get the same recursive
formulas for both cases.\\

Define $$P_{\grn}(q,t,s)=\sumlim_{\pi \in \grn}q^{{\rm
exc}(\pi)}t^{{\rm fix}(\pi)}s^{{\rm cyc}(\pi)}.$$

Concerning $\grn$, we prove the following two main results:

\bthm \label{theorem1}

$$P_{\grn}^{{\rm Abs}}(q,1,-1)=P_{\grn}^{{\rm Clr}}(q,1,-1)=(q^r-1)P_{G_{r,n-1}}(q,1,-1).$$
Hence,
$$P^{{\rm Abs}}_{\grn}(q,1,-1)=P^{\rm Clr}_{\grn}(q,1,-1)=-\frac{(q^r-1)^n}{q-1}.$$ \ethm

\bthm \label{theorem2}
$$P_{\grn}^{{\rm Abs}}(q,0,-1)=P_{\grn}^{{\rm Clr}}(q,0,-1)=[r]_q( P_{G_{r,n-1}}(q,0,-1) -q^{n-1} [r]_q^{n-1}).$$
Hence,
$$P^{\rm Abs}_{\grn} (q,0,-1)=P^{\rm Clr}_{\grn} (q,0,-1)= -q [r]_q^n [n-1]_q.$$
\ethm

One can easily check that the formulas appeared in Theorem
\ref{theorem1} and Theorem \ref{theorem2} indeed
 generalize the formulas of Ksavrelof and Zeng (for $r=1$).\\

We apply our techniques also to obtain permutations statistics on
the group of even signed permutations, $D_n$, also known as the
{\it Coxeter group of type $D$.} We get the following results:

\bthm \label{theorem3}
$$P_{D_n}^{{\rm Clr}}(q,1,-1)=(q^2-1)P_{D_{n-1}}^{{\rm Clr}}(q,1,-1).$$
Hence,
$$P^{\rm Clr}_{D_n}(q,1,-1)=(1-q^2)^{n-1}.$$
\ethm

\bthm \label{theorem4}
$$P_{D_n}^{{\rm Abs}}(q,1,-1)=-\frac{1}{2}(q-1)^{n-1}((1+q)^n+(1-q)^n).$$
\ethm

This paper is organized as follows. In Section \ref{pre}, we
recall some properties of $\grn$. In Section \ref{stat} we define
the new statistics on $\grn$. Sections \ref{proof theorem 1} and
\ref{proof theorem 2} deal with the proofs of Theorems
\ref{theorem1} and \ref{theorem2} respectively. Section \ref{D_n}
includes the proofs of Theorems \ref{theorem3} and \ref{theorem4}.

\section{Preliminaries}\label{pre}

\subsection{Notations}\label{notations}

 For $n \in \mathbb{N}$, let $[n] \eqdef \{ 1,2, \ldots , n \} $
(where $[0] \eqdef \emptyset $).

Also, let:
$$[n]_q:=\frac{1-q^n}{1-q}=1+q+\cdots +q^{n-1}$$
(so $[0]_q=0$), and
$$[n]_q!=[n]_q \cdot [n-1]_q \cdots [1]_q.$$

\subsection{The group of colored permutations}\label{grn}

\bde
Let $r$ and $n$ be positive integers. {\it The group of colored
permutations of $n$ digits with $r$ colors} is the wreath product
$$\grn=\mathbb{Z}_r \wr S_n=\mathbb{Z}_r^n \rtimes S_n,$$ consisting of all the pairs $(z,\tau)$ where
$z$ is an $n$-tuple of integers between $0$ and $r-1$ and $\tau
\in S_n$. The multiplication is defined by the following rule: For
$z=(z_1,...,z_n)$ and $z'=(z'_1,...,z'_n)$
$$(z,\tau) \cdot (z',\tau')=((z_1+z'_{\tau(1)},...,z_n+z'_{\tau(n)}),\tau \circ \tau')$$ (here $+$ is
taken modulo $r$).
\ede

Here are some conventions we use along this paper:
For an  element  $\pi=(z,\tau) \in \grn$ with $z=(z_1,...,z_n)$ we
write $z_i(\pi)=z_i$.
For $\pi=(z,\tau)$, we denote $|\pi|=(0,\tau), (0 \in \mathbb{Z}_r^n)$.
An element $(z,\tau)=((1,0,3,2),(2,1,4,3)) \in G_{3,4}$
will be written as $(\bar{2} 1 \bar{\bar{\bar{4}}} \bar{\bar{3}})$.

A much more natural way to present $\grn$ is the following:
Consider the alphabet $\Sigma=\{1,\dots,n,\bar{1},\dots,\bar{n},\dots,
1^{[r-1]},\dots,n^{[r-1]} \}$ as the set $[n]$ colored by the colors
$0,\dots,r-1$. Then, an element of $\grn$ is a {\it colored
permutation}, i.e. a bijection $\pi: \Sigma \rightarrow \Sigma$ such that
$\pi(\bar{i})={\overline{\pi (i)}}$.

In particular, $G_{1,n}=C_1 \wr S_n$ is the symmetric group $S_n$,
while $G_{2,n}=C_2\wr S_n$ is the group of signed permutations
$B_n$, also known as the {\it hyperoctahedral group}, or the {\it
classical Coxeter group of type B}. We also define here the
following normal subgroup of $B_n$ of index $2$, called the {\it
even signed permutation group} or the {\it Coxeter group of type
$D$}:
$$D_n=\{\pi \in B_n \mid \sumlim_{i=1}^n{z_i(\pi)}\equiv 0 \pmod 2\}.$$

\section{Statistics on $\grn$}\label{stat}

Given any ordered alphabet $\Sigma'$, we recall the definition of
the {\it excedance set} of a permutation $\pi$ on $\Sigma'$ :
$${\rm Exc}(\pi)=\{i \in \Sigma' \mid \pi(i)>i\}$$ and the {\it excedance
number} is defined to be ${\rm exc}(\pi)=|{ \rm Exc}(\pi)|$.\\

We start by defining two orders on the set:
$$\Sigma=\{1,\dots,n,\bar{1},\dots,\bar{n},\dots,
1^{[r-1]},\dots,n^{[r-1]} \}.$$

\bde The {\it absolute order} on $\Sigma$ is defined to be:
$$1^{[r-1]} < \cdots < \bar{1} < 1 < 2^{[r-1]} < \cdots < \bar{2} <2 < \cdots < n^{[r-1]} < \cdots < \bar{n} < n.$$
The {\it color order} on $\Sigma$ is defined to be:
$$1^{[r-1]} < \cdots < n^{[r-1]} < 1^{[r-2]} < 2^{[r-2]} < \cdots < n^{[r-2]} < \cdots < 1 < \cdots < n.$$
\ede

\begin{exa}

Given the color order: $$\bar{\bar{1}} < \bar{\bar{2}}
<\bar{\bar{3}} < \bar{1} < \bar{2} <\bar{3} < 1 < 2 < 3,$$ we
write $\sigma=(3\bar{1}\bar{\bar{2}}) \in G_{3,3}$ in an extended
form:

$$\begin{pmatrix} \bar{\bar{1}} & \bar{\bar{2}} & \bar{\bar{3}} &
\bar{1} & \bar{2}& \bar{3} & 1 & 2 & 3\\
\bar{\bar{3}} & 1 & \bar{2} & \bar{3} & \bar{\bar{1}}  &  2 & 3 &
\bar{1} & \bar{\bar{2}}
\end{pmatrix}$$
and calculate:
${\rm Exc}(\sigma)=\{\bar{\bar{1}},\bar{\bar{2}},\bar{\bar{3}},\bar{1},\bar{3},1\}$
and ${\rm exc}(\sigma)=6$.
\end{exa}

Before defining the excedance number with respect to both orders, we have to introduce some notions.

Let $\sigma \in \grn$. We define:
$${\rm csum}(\sigma) = \sumlim_{i=1}^n z_i(\sigma),$$

$${\rm Exc}_A(\sigma) = \{ i \in [n-1] \ | \ \sigma(i) > i \},$$
where the comparison is with respect to the color order.
$${\rm exc}_A(\sigma) = |{\rm Exc}_A(\sigma)|.$$

\begin{exa}
Take $\sigma=(\bar{1}\bar{\bar{3}}4\bar{2}) \in G_{3,4}$. Then ${\rm csum}(\sigma)=4$,
${\rm Exc_A}(\sigma)=\{3\}$ and hence ${\rm exc}_A(\sigma)=1$.
\end{exa}

Let $\sigma \in \grn$. Recall that for $\sigma=(z,\tau)\in \grn$,
$|\sigma|$ is the permutation of $[n]$ satisfying
$|\sigma|(i)=\tau(i)$. For example, if
$\sigma=(\bar{2}\bar{\bar{3}}1\bar{4})$ then $|\sigma|=(2314)$.

Now we can define the excedance numbers for $\grn$.

\bde
Define:
$${\rm exc}^{{\rm Abs}}(\sigma)={\rm exc}(|\sigma|)+{\rm csum}(\sigma),$$
$${\rm exc}^{{\rm Clr}}(\sigma)=r \cdot {\rm exc}_A(\sigma)+{\rm csum}(\sigma).$$
\ede

The parameters
${\rm exc}^{\rm Abs}$ and ${\rm exc}^{\rm Clr}$ are indeed
different:  for $\sigma=(21) \in G_{r,2}$, $(r>1)$  one has
${\rm exc}^{\rm Abs}(\sigma)=1$ but
${\rm exc}^{\rm Clr}(\sigma)=r$.

One can view ${\rm exc}^{{\rm Clr}}(\sigma)$ in a different way:

\blem Let $\sigma \in \grn$. Consider the set $\Sigma$ ordered
by the color order. Then $${\rm exc}(\sigma)={\rm exc}^{{\rm
Clr}}(\sigma).$$ \elem

\begin{proof}
Let $i \in [n]$. We divide our proof into two cases:
$z_i(\sigma)=0$ and $z_i(\sigma) \neq 0$.

If $z_i(\sigma)=0$, then $i \in {\rm Exc}_A (\sigma)$ if and only
if $\sigma(i)>i$. In this case, we have $\sigma(i^{[j]})> i^{[j]}$
for every color $1 \leq j \leq r-1$. Hence, we have $\{ i,i^{[1]},
\cdots, i^{[r-1]}\} \subseteq {\rm Exc} (\sigma)$. Hence, each $i
\in {\rm Exc}_A (\sigma)$ contributes $r$ excedances to ${\rm exc}
(\sigma)$.

On the other hand, if $z_i(\sigma) = k \neq 0$, we have that $i
\notin {\rm Exc}(\sigma)$. By definition, we have
$\sigma(i^{[j]})= |\sigma (i)|^{[(j+k) \pmod r]}$ for all $j$.
Thus, for $0 \leq j \leq r-k-1$, $i^{[j]} \notin {\rm Exc}
(\sigma)$, and for the $k$ indices $r-k \leq j \leq r-1$, $i^{[j]}
\in {\rm Exc} (\sigma)$.

Consequently, we have: $${\rm exc}(\sigma)={\rm exc}^{{\rm
Clr}}(\sigma).$$

\end{proof}

Recall that any permutation of $S_n$ can be decomposed into a
product of disjoint cycles. This notion can be easily generalized
to the group $\grn$ as follows. Given any $\pi \in \grn$ we define
the {\it cycle number} of $\pi=(z,\tau)$ to be the number of
cycles in $\tau$.

We say that $i \in [n]$ is an {\it absolute fixed point} of
$\sigma \in \grn$ if $|\sigma(i)|=i$.

\section{Proof of Theorem \ref{theorem1}}\label{proof theorem 1}

In this section we prove Theorem \ref{theorem1}. The idea of proving this type of identities is
constructing a subset $S$ of $\grn$ whose contribution to the generating function is exactly
the right side of the identity and a killing involution on $\grn  - S$, i.e., an involution
on $\grn - S$ which preserves the number of excedances but changes the sign of every element of
$\grn -S$ and hence shows that $\grn -S$ contributes nothing to the generating function.

\subsection{Proof for the absolute order}\label{ProofAbs1}
We divide $\grn$ into $2r+1$ disjoint subsets as follows:

$$K_{r,n}=\{\sigma \in \grn \mid |\sigma(n)|\neq n , |\sigma(n-1)| \neq n\}.$$

$$T^{i}_{r,n}=\{\sigma \in \grn \mid \sigma(n)=n^{[i]}\}, \qquad
(0 \leq i \leq r-1).$$

$$R^{i}_{r,n}=\{\sigma \in \grn \mid \sigma(n-1)=n^{[i]}\}, \qquad
(0 \leq i \leq r-1).$$

We first construct a killing involution on the set $K_{r,n}$. Let
$\sigma \in K_{r,n}$. Define $\varphi: K_{r,n} \to K_{r,n}$ by:
$$\sigma ' = \varphi(\sigma)=(\sigma(n-1),\sigma(n)) \sigma.$$
Note that $\varphi$ exchanges $\sigma(n-1)$ with $\sigma(n)$. It is obvious
that $\varphi$ is indeed an involution.

We will show that $\rm{exc}^{{\rm Abs}}(\sigma)=\rm{exc}^{{\rm
Abs}}(\sigma ')$. First, for $i<n-1$, it is clear that $i \in
{\rm Exc}(|\sigma|)$ if and only if $i \in {\rm Exc}(|\sigma'|)$. Now, as
$\sigma(n-1) \neq n$, $n-1 \notin {\rm Exc}(|\sigma|)$ and thus $n
\notin {\rm Exc}(|\sigma'|).$ Finally, $|\sigma(n)|\neq n$ implies
that $n-1 \notin {\rm Exc}(|\sigma'|)$ and thus
${\rm exc}^{{\rm Abs}}(\sigma)={\rm exc}^{{\rm Abs}}(\sigma').$

On the other hand, $\rm{cyc}(\sigma)$ and $\rm{cyc}(\sigma')$ have
different parities due to a multiplication by a transposition.
Hence, $\varphi$ is indeed a killing involution on $K_{r,n}$.\\

We turn now to the sets $T^i_{r,n} \quad (0 \leq i \leq r-1)$ . Note
that there is a natural bijection between $T^i_{r,n}$ and
$G_{r,n-1}$ defined by ignoring the last digit. Let $\sigma \in
T^{i}_{r,n}$. Denote the image of $\sigma \in T^{i}_{r,n}$ under
this bijection by $\sigma'$. Since $n \not\in \rm{Exc}(|\sigma|)$,
we have $\rm{exc}(|\sigma|)=\rm{exc}(|\sigma'|)$. Now,
$\rm{csum}(\sigma')=\rm{csum}(\sigma)-i$, since $z_n(\sigma)=i$
and hence we have:
$$\rm{exc}^{{\rm
Abs}}(\sigma)-i=\rm{exc}^{{\rm
Abs}}(\sigma').$$ Finally, since $n$ is an
absolute fixed point of $\sigma$, $\rm{cyc}(\sigma ')=
\rm{cyc}(\sigma)-1$ and we get that the total contribution of
$T^i_{r,n}$ is:
$$P_{T_{r,n}^i}^{{\rm Abs}}(q,1,-1) =-q^i P_{G_{r,n-1}}^{{\rm Abs}}(q,1,-1)$$
for $0 \leq i \leq r-1$.\\

Now, we treat the sets $R^i_{r,n} \quad (0 \leq i \leq r-1)$. There is a
bijection between $R^i_{r,n}$ and $T^i_{r,n}$ using the same
function $\varphi$ we used above. Let $\sigma \in R_{r,n}^i$.
Define $\varphi: R^i_{r,n} \to T^i_{r,n}$ by:
$$\sigma ' = \varphi(\sigma)=(\sigma(n-1),\sigma(n)) \sigma.$$

 In $\sigma$, we have that $n-1 \in
\rm{Exc}(|\sigma|)$ (since $|\sigma(n-1)|=n$) and $n \not\in
\rm{Exc}(|\sigma|)$, but in $\sigma'$, $n-1,n \not\in
\rm{Exc}(|\sigma'|)$. Hence, $\rm{exc}
(|\sigma|)-1=\rm{exc}(|\sigma'|)$. We also have that
$\rm{csum}(\sigma) =\rm{csum}(\sigma')$ (since
$z_{n-1}(\sigma)+z_n(\sigma) = z_{n-1}(\sigma')+z_n(\sigma')$).
Hence, we have that
$$\rm{exc}^{{\rm Abs}}(\sigma)-1=\rm{exc}^{{\rm Abs}}(\sigma ').$$

As before, the number of cycles changes its parity due to the
multiplication by a transposition, and thus:
$(-1)^{\rm{cyc}(\sigma)} = -(-1)^{\rm{cyc}(\sigma')}$.\\

Hence, the total contribution of the elements in $R_{r,n}^i$ is
$$P_{R_{r,n}^i}^{{\rm Abs}}(q,1,-1)=q^{i+1} \cdot P_{G_{r,n-1}}^{{\rm Abs}}(q,1,-1)$$
for $0 \leq i \leq r-1$.

Now, if we sum up all the parts, we get:

\begin{eqnarray*}
P_{\grn}^{{\rm Abs}} (q,1,-1) &=& P_{K_{r,n}}^{{\rm Abs}}(q,1,-1)
+\sumlim_{i=0}^{r-1}{P_{T_{r,n}^i}^{{\rm Abs}}(q,1,-1)}
+\sumlim_{i=0}^{r-1}{P_{R_{r,n}^i}}(q,1,-1)\\
&=&  \sum_{i=0}^{r-1} (-q^i P_{G_{r,n-1}}^{{\rm Abs}}(q,1,-1)) +
\sum_{i=0}^{r-1} q^{i+1} P_{G_{r,n-1}}^{{\rm Abs}}(q,1,-1) \\&=&
(q^r -1) P_{G_{r,n-1}}^{{\rm Abs}}(q,1,-1)
\end{eqnarray*}
as claimed.\\

Now, for $n=1$, $G_{r,1}$ is the cyclic group of order $r$ and
thus
$$P^{\rm{Abs}}_{G_{r,1}}(q,1,-1)=-(1+q+\cdots +q^{r-1})=-\frac{q^r-1}{q-1},$$ so we have
$$P_{\grn}^{{\rm Abs}}(q,1,-1)=-\frac{(q^r-1)^n}{{q-1}}.$$

\subsection{Proof for the color order}

As in the previous proof, we divide $\grn$ into the same $2r+1$
disjoint subsets $K_{r,n}$, $T^{i}_{r,n} \quad (0 \leq i \leq
r-1)$ and $R^{i}_{r,n} \quad (0 \leq i \leq r-1)$ used there.

As before, we first construct a killing involution on the set $K_{r,n}$. Let
$\sigma \in K_{r,n}$. As before, define $\varphi: K_{r,n} \to K_{r,n}$ by :
$$\sigma ' = \varphi(\sigma)=(\sigma(n-1),\sigma(n)) \sigma.$$
The proof that $\varphi$ is a killing involution is similar to the
one we presented in Section \ref{ProofAbs1}.\\

We turn now to the sets $T^i_{r,n}$. We use again the bijection
between $T^i_{r,n}$ and $G_{r,n-1}$ defined by ignoring the last
digit. Let $\sigma \in T_{r,n}$. As in the previous proof, we have:
$$\rm{exc}^{{\rm Clr}}(\sigma)-i=\rm{exc}^{{\rm Clr}}(\sigma').$$
Now, since $n$ is an absolute fixed point of $\sigma$,
$\rm{cyc}(\sigma ')= \rm{cyc}(\sigma)-1$.

To summarize, we get that the total contribution of elements
in $T^i_{r,n}$ is:
$$P_{T_{r,n}^i}^{\rm Clr}=-q^i P_{G_{r,n-1}}^{{\rm Clr}}(q,1,-1)$$
for $0 \leq i \leq r-1$.\\

Finally, we treat the sets $R^i_{r,n}$. Let $\sigma \in R^i_{r,n}$.
Recall the bijection $\varphi:R^i_{r,n} \rightarrow T_{r,n}^i$
defined in Section \ref{ProofAbs1} by:
$$\sigma ' = \varphi(\sigma)=(\sigma(n-1),\sigma(n)) \sigma.$$

When we compute the change in the excedance, we split our
treatment into two cases: $i=0$ and $i>0$. For the case $i=0$, we
get $\rm{exc}^{{\rm Clr}}(\sigma)-r=\rm{exc}^{{\rm Clr}}(\sigma
')$. For the case $i>0$, we show that $\rm{exc}^{{\rm
Clr}}(\sigma)=\rm{exc}^{{\rm Clr}}(\sigma ')$.

We start with case $i=0$. Note that $n-1 \in \rm{Exc}_A(\sigma)$
(since $\sigma(n-1)=n$) and $n \not\in \rm{Exc}_A(\sigma)$. On the
other hand, in $\sigma'$, $n-1,n \not\in \rm{Exc}_A(\sigma')$. Hence,
$\rm{exc}_A (\sigma)-1=\rm{exc}_A(\sigma')$.

Now, for the case $i>0:$ $n-1,n \not\in \rm{Exc}_A(\sigma)$ (since
$\sigma(n-1)=n^{[i]}$ is not an excedance with respect to the
color order). We also have: $n-1,n \not\in \rm{Exc}_A(\sigma')$ and
thus $\rm{Exc}_A(\sigma)=\rm{Exc}_A(\sigma ')$ for $\sigma \in
R^i_{r,n}$ where $i>0$.

In both cases, we have that $\rm{csum}(\sigma)
=\rm{csum}(\sigma')$. Hence, we have that $\rm{exc}^{{\rm
Clr}}(\sigma)-r=\rm{exc}^{{\rm Clr}}(\sigma ')$ for $i=0$ and
$\rm{exc}^{{\rm Clr}}(\sigma)=\rm{exc}^{{\rm Clr}}(\sigma ')$ for
$i>0$.

As before, the number of cycles changes its parity due to the multiplication
by a transposition, and hence:
$(-1)^{\rm{cyc}(\sigma)} = -(-1)^{\rm{cyc}(\sigma')}$.

Hence, the total contribution of elements in $R^i_{r,n}$ is:
$$q^r P_{G_{r,n-1}}^{{\rm Clr}}(q,1,-1)$$
for $i=0$, and
$$q^i P_{G_{r,n-1}}^{{\rm Clr}}(q,1,-1)$$
for $i>0$.\\

Now, if we sum up all the parts, we get:

$$P_{\grn}^{{\rm Clr}} (q,1,-1)=$$ $$\sum_{i=0}^{r-1} (-q^i
P_{G_{r,n-1}}^{{\rm Clr}}(q,1,-1)) + q^r P_{G_{r,n-1}}^{{\rm
Clr}}(q,1,-1)  +  \sum_{i=1}^{r-1} q^i P_{G_{r,n-1}}^{{\rm
Clr}}(q,1,-1)=$$  $$(q^r -1) P_{G_{r,n-1}}^{{\rm Clr}}(q,1,-1)$$

as needed.\\

Now, for $n=1$, $G_{r,1}$ is the cyclic group of order $r$ and
thus
$$P^{\rm{Clr}}_{G_{r,1}}(q,1,-1)=-(1+q+\cdots +q^{r-1})=-\frac{q^r-1}{q-1},$$ so we have
$$P^{\rm{Clr}}_{\grn}(q,1,-1)=-\frac{(q^r-1)^n}{{q-1}}.$$

\section{Derangements in $\grn$ and the proof of Theorem \ref{theorem2}}\label{proof theorem 2}

We start with the definition of a derangement.

\bde An element $\sigma \in \grn$ is called a {\it derangement} if
it has no absolute fixed points, i.e. $|\sigma(i)|\neq i$ for
every $i \in [n]$. Denote by $D_{r,n}$ the set of all derangements
in $\grn$. \ede

In this section, we prove Theorem \ref{theorem2}. As in the
previous section, we prove Theorem \ref{theorem2} for both
orders.

\subsection{Proof for the absolute order}

We divide $D_{r,n}$ into $r+2$ disjoint subsets in the following
way:

$$A^i_{r,n}=\{\sigma \in D_{r,n} \mid \sigma(2)=1^{[i]},|\sigma(1)| \neq 2\}, \qquad i=0,\dots ,r-1.$$

$$T_{r,n}=\{\sigma \in D_{r,n} \mid |\sigma|=(234 \cdots n1)\}.$$

$$\hat{D}_{r,n}=D_{r,n} - (\bigcup_{i=0}^{r-1}{A^i_{r,n}} \cup T_{r,n}).$$

We start by constructing a killing involution $\varphi$ on
$\hat{D}_{r,n}$. Given any $\sigma \in \hat{D}_{r,n}$, let $i$ be the first number
such that $|\sigma(i)| \neq i+1$. Define
$$\sigma'=\varphi(\sigma)=(\sigma(i),\sigma(i+1)) \sigma.$$
For example, if $\sigma=(\bar{3}4\bar{1}\bar{5}\bar{\bar{2}})$ then
$\sigma'=(4\bar{3}\bar{1}\bar{5}\bar{\bar{{2}}})$.

It is easy to see that $\varphi$ is a well-defined involution on
$\hat{D}_{r,n}$. We proceed to prove that ${\rm exc}^{{\rm
Abs}}(\sigma)={\rm exc}^{{\rm Abs}}(\sigma')$. Indeed, ${\rm
csum}(\sigma)={\rm csum}(\sigma')$.

Let $i$ be the first number such that $|\sigma(i)| \neq i+1$ so
that in the pass from $\sigma$ to $\sigma'$ we exchange
$\sigma(i)$ with $\sigma(i+1)$. For every $j \neq i,i+1$, clearly
$j \in {\rm Exc}(|\sigma|)$ if and only if $j \in {\rm
Exc}(|\sigma'|)$. Since $\sigma \in D_{r,n}$, $|\sigma(i)| \neq
i+1$ and $|\sigma(j)|=j+1$ for $j<i$, we have that
$|\sigma(i)|,|\sigma(i+1)| \in \{1,i+2,\cdots,n\}$. Thus,
exchanging $\sigma(i)$ with $\sigma(i+1)$ does not change ${\rm
exc}(|\sigma|)$.\\

Note also that the parity of ${\rm cyc}(\sigma')$ is opposite to
the parity of ${\rm cyc}(\sigma)$ due to the multiplication by
a transposition. Hence, we have proven that
$\varphi$ is a killing involution.

Now, let us calculate the contribution of each set in our
decomposition to $P_{D_{r,n}}^{{\rm Abs}}(q,0,-1)$. As we have
shown, $\hat{D}_{r,n}$ contributes nothing. Define a bijection
$$\psi: A^{i}_{r,n} \to D_{r,n-1}$$ by: $\psi(\sigma)=\sigma'$
where $\sigma'(1)=(|\sigma(1)|-1)^{z_1(\sigma)}$ and for $i>1$,
$\sigma'(i)=(|\sigma(i+1)|-1)^{z_{i+1}(\sigma)}$. For example, if
$\sigma=(3\bar{1}4\bar{2})$, then $\sigma '=(23\bar{1})$. It is
easy to see that ${\rm exc}(|\sigma|)={\rm exc}{(|\sigma'|)}$. On
the other hand, ${\rm csum}(\sigma')={\rm csum}(\sigma)-i$ and
${\rm cyc}(\sigma)={\rm cyc}(\sigma')$ and thus the contribution
of $A^{i}_{r,n}$ to $P_{D_{r,n}}^{{\rm Abs}}(q,0,-1)$ is $q^i
P_{D_{r,n-1}}^{{\rm Abs}}(q,0,-1)$ for $1 \leq i \leq r-1$.\\

Finally, we treat the set $T_{r,n}$. For every $\sigma \in
T_{r,n}$ we have: ${\rm exc}(|\sigma|)=n-1$ and ${\rm
cyc}(\sigma)=1$. Concerning ${\rm csum}(\sigma)$ we have:
$$\sumlim_{\sigma \in T_{r,n}}{q^{{\rm
csum}(\sigma)}}=(1+q+\cdots q^{r-1})^n.$$

To summarize, we get

\begin{eqnarray*}
P_{D_{r,n}}^{{\rm Abs}} (q,0,-1) & = & \sum_{i=0}^{r-1} (q^i
P_{D_{r,n-1}}^{{\rm Abs}}(q,0,-1)) +q^{n-1}(1+q+\cdots q^{r-1})^n   \\
& = & (1+q+\cdots q^{r-1})P_{D_{r,n-1}}^{{\rm
Abs}}(q,0,-1)+q^{n-1}(1+q+ \cdots q^{r-1})^n  \\
& = & [r]_q(P_{D_{r,n-1}}^{{\rm Abs}}(q,0,-1)-q^{n-1}[r]_q^{n-1})
\end{eqnarray*} as needed.\\

Now, for $n=2$, we have:
$$P^{\rm Abs}_{D_{r,2}}(q,0,-1) =-q[r]_q^2$$
and thus
$$P^{\rm Abs}_{D_{r,n}}(q,0,-1)=-q[r]_q^n[n-1]_q.$$

\subsection{Proof for the color order}

We use the same decomposition of $D_{r,n}$ as before. The killing
involution will be also the same, due to the following
observation: one can replace ${\rm Exc}(|\sigma|)$ in the previous
proof by ${\rm Exc}_A(\sigma)$, and the argument still holds. Note
that $i \notin {\rm Exc}_A(\sigma)$ if $z_i(\sigma) \neq 0$.

Now, let us calculate the contribution of each set in our
decomposition to $P_{D_{r,n}}^{{\rm Clr}}(q,0,-1)$. As we have
shown, $\hat{D}_{r,n}$ contributes nothing.\\

As before, define a bijection from $A^{i}_{r,n}$ to $D_{r,n-1}$ by
$\sigma \mapsto \sigma'$ where $\sigma'(1)=(|\sigma(1)|-1)^{z_1(\sigma)}$ and for $i>1$,
$\sigma'(i)=(|\sigma(i+1)|-1)^{z_{i+1}(\sigma)}$. As before, it is
easy to see that ${\rm exc}_A(\sigma)={\rm exc}_A(\sigma')$ (since
$|\sigma(2)|=1$). On the other hand, ${\rm csum}(\sigma')={\rm
csum}(\sigma)-i$ and ${\rm cyc}(\sigma)={\rm cyc}(\sigma')$, and
thus the contribution of $A^{i}_{r,n}$ to $P_{D_{r,n}}^{{\rm
Clr}}(q,0,-1)$ is $q^i P_{D_{r,n-1}}^{{\rm Clr}}(q,0,-1)$ for $0
\leq i \leq r-1$.\\

Now, we treat the set $T_{r,n}$. Let $\sigma \in T_{r,n}$. Observe
that for $i<n$, $i \in {\rm Exc}_A(\sigma)$ if and only if
$z_i(\sigma)=0$. In this case, the place $i$ contributes $r$ to
${\rm exc}^{{\rm Clr}}(\sigma)$. Hence, it will be natural to
construct a bijection between $T_{r,n}$ and the following subset
of $G_{r+1,n}$:
$$W = \{ \sigma \in G_{r+1,n} \ | \ |\sigma|=(23 \cdots n1), z_i(\sigma)\neq 0\ {\rm for}\ 1\leq i \leq n-1, z_n(\sigma) \neq r \}.$$
The bijection is defined by the rule
$\sigma \longmapsto \sigma '$ where:
$$
\sigma'(i) =
\left\{
\begin{array}{cc}
  \sigma(i) & i=n\ {\rm or}\ z_i(\sigma)\neq 0 \\
  |\sigma(i)|^{[r]} & {\rm otherwise.}
\end{array}
\right.
$$
Note that ${\rm exc}^{{\rm Clr}}(\sigma)= {\rm csum}(\sigma')$.

Now we compute:
\begin{eqnarray*}
\sumlim_{\sigma \in T_{r,n}}{q^{{\rm exc}^{{\rm Clr}}(\sigma)}(-1)^{{\rm cyc}(\sigma)}} &=&
\sumlim_{\sigma' \in W}{q^{{\rm csum}(\sigma')}(-1)^{{\rm cyc}(\sigma')}} \\ &=&
-\sumlim_{\sigma' \in W}{q^{{\rm csum}(\sigma ')}}\\&=&
-(q+q^2+\cdots +q^r)^{n-1}(1+q+\cdots+q^{r-1}).
\end{eqnarray*}

To summarize, we get

\begin{eqnarray*}
P_{D_{r,n}}^{{\rm Clr}} (q,0,-1) & = & \sum_{i=0}^{r-1} (q^i
P_{D_{r,n-1}}^{{\rm Clr}}(q,0,-1))-(q+q^2+\cdots +q^r)^{n-1}(1+q+\cdots+q^{r-1}) \\ & = &
(1+q+\cdots q^{r-1})(P_{D_{r,n-1}}^{{\rm Clr}}(q,0,-1)+(q+ \cdots q^r)^{n-1}) \\ & = &
[r]_q(P_{D_{r,n-1}}^{{\rm Clr}}(q,0,-1)-q^{n-1}[r]_q^{n-1}).
\end{eqnarray*}
Now, as before, for $n=2$, we have:
$$P^{\rm Clr}_{D_{r,2}}(q,0,-1) =-q[r]_q^2$$
and thus

$$P^{\rm Clr}_{D_{r,n}}(q,0,-1)=-q[r]_q^n[n-1]_q.$$

\section{Statistics on the group of even signed permutations} \label{D_n}
In this section we deal with the Coxeter group of type $D$, namely
the group of even signed permutations. We recall its
definition:
$$D_n=\{\pi \in B_n \mid \sumlim_{i=1}^n{z_i(\pi)}\equiv 0 \pmod 2\}.$$

Unlike the case of the groups $\grn$, in $D_n$ the distribution of the excedance numbers with respect
to the color order is different from the distribution with respect to the absolute order. We start with the color order.

\subsection{Proof of Theorem \ref{theorem3}}

We divide $D_n$ into $5$ subsets:
$$K_n=\{\sigma \in D_n \mid |\sigma(n)|\neq n,|\sigma(n-1)|\neq
n\}.$$

$$T_n^0=\{\sigma \in D_n \mid \sigma (n)=n\}.$$
$$T_n^1=\{\sigma \in D_n \mid \sigma (n)=\bar{n}\}.$$
$$R_n^0=\{\sigma \in D_n \mid \sigma (n-1)=n\}.$$
$$R_n^1=\{\sigma \in D_n \mid \sigma (n-1)=\bar{n}\}.$$

We denote:
$$a_n=P_{D_n}^{\rm Clr}(q,1,-1),$$
$$b_n=P_{D_n^c}^{\rm Clr}(q,1,-1),$$
 where $D_n^c$ is the complement of $D_n$ in $B_n$.

Define $\varphi: K_n \to K_n$ by
$$\sigma ' = \varphi(\sigma)=(\sigma(n-1),\sigma(n)) \sigma.$$
Note that $\varphi$ exchanges $\sigma(n-1)$ with $\sigma(n)$. It
is easy to see that $\varphi$ is a killing involution on $K_n$.

We turn now to the set $T_n^0$. Note that there is a natural
bijection between $T_n^0$ and $D_{n-1}$, defined by ignoring the
last digit. Let $\sigma \in T_n^0$. Denote the image of $\sigma$ under this
bijection by $\sigma'$. Note that ${\rm
csum}(\sigma')={\rm csum}(\sigma)$, ${\rm Exc_A}(\sigma')={\rm
Exc_A}(\sigma)$ and ${\rm Exc}^{{\rm Clr}}(\sigma')={\rm
Exc}^{{\rm Clr}}(\sigma)$. On the other hand, ${\rm
cyc}(\sigma')={\rm cyc}(\sigma)-1$ and thus the restriction of
$a_n$ to $T_n^0$ is just $-a_{n-1}$.

For the contribution of the set $T_n^1$, note that the function
$\varphi$ defined above gives us a bijection between $T_n^1$ and
$D_{n-1}^c$. In this case, ${\rm csum}(\sigma')={\rm
csum}(\sigma)-1$, ${\rm exc_A}(\sigma')={\rm exc_A}(\sigma)$ and
${\rm exc}^{{\rm Clr}}(\sigma')={\rm exc}^{{\rm Clr}}(\sigma)$. On
the other hand, ${\rm cyc}(\sigma')={\rm cyc}(\sigma)-1$ as
before. Hence, the restriction of $a_n$ to $T_n^1$ is $-q
b_{n-1}$.

Now, for the set $R_n^0$, we have the following bijection between
$R_n^0$ and $D_{n-1}$: for $\sigma \in R_n^0$, exchange the last
two digits, and then ignore the last digit. If we denote the image
of $\sigma$ by $\sigma'$, we have ${\rm csum}(\sigma')={\rm
csum}(\sigma)$, ${\rm exc_A}(\sigma')={\rm exc_A}(\sigma)-1$,
${\rm exc}^{{\rm Clr}}(\sigma')={\rm exc}^{{\rm Clr}}(\sigma)-2$
and ${\rm cyc}(\sigma') \equiv {\rm cyc}(\sigma) \pmod 2$. Hence,
the restriction of $a_n$ to $R_n^0$ is $q^2 a_{n-1}$.

For the set $R_n^1$, we have a bijection between $R_n^1$ and
$D_{n-1}^c$: for $\sigma \in R_n^1$, exchange the last two digits,
and then ignore the last digit. Denoting the image of $\sigma$ by
$\sigma'$, we have ${\rm csum}(\sigma')={\rm csum}(\sigma)-1$,
${\rm exc_A}(\sigma')={\rm exc_A}(\sigma)$, and hence ${\rm
exc}^{{\rm Clr}}(\sigma')={\rm exc}^{{\rm Clr}}(\sigma)-1$. Also, we have ${\rm
cyc}(\sigma') \equiv {\rm cyc}(\sigma) \pmod 2$. Hence, the
restriction of $a_n$ to $R_n^1$ is $q b_{n-1}$.

We summarize all the contributions over all the four subsets, and
we have:
$$a_n=-a_{n-1}-q b_{n-1}+q^2 a_{n-1} +q b_{n-1}=(q^2-1)a_{n-1}$$

For computing $a_1$, note that $D_1 =\{1\}$ and thus $a_1=-1$.

Therefore, we have:
$$P_{D_n}^{\rm Clr}(q,1,-1)=a_n=-(q^2-1)^{n-1},$$
and we are done.

\subsection{Proof of Theorem \ref{theorem4}}
In this subsection we present the proof of Theorem \ref{theorem4}
which computes $P_{D_n}^{{\rm Abs}}(q,1,-1)$.\\

We start by dividing $D_n$ into $5$ subsets just as was shown in the proof of Theorem \ref{theorem3} and
define as before:

$$a_n=P_{D_n}^{\rm Abs}(q,1,-1),$$
$$b_n=P_{D_n^c}^{\rm Abs}(q,1,-1),$$
 where $D_n^c$ is the complement of $D_n$ in $B_n$.

It is easy to check that the sets $T_n^0$ and $T_n^1$ give the
same contributions as before, so we turn to the set $R_n^0$. By
using the bijection between $R_n^0$ and $D_{n-1}$ defined above
which exchanges the last two digits, and then ignores the last
digit, we have: ${\rm csum}(\sigma')={\rm csum}(\sigma)$, ${\rm
exc}(|\sigma'|)={\rm exc}(|\sigma|)-1$. Thus, ${\rm exc}^{{\rm
Abs}}(\sigma ')={\rm exc}^{{\rm Abs}}(\sigma)-1$. Also, ${\rm
cyc}(\sigma') \equiv {\rm cyc}(\sigma) \pmod 2$. Hence, the
restriction of $a_n$ to $R_n^0$ is $q a_{n-1}$.

For the set $R_n^1$, we use the same bijection, now between
$R_n^1$ and $D_{n-1}^c$ to get: ${\rm csum}(\sigma')={\rm
csum}(\sigma)-1$, ${\rm exc}(|\sigma'|)={\rm exc}(|\sigma|)-1$.
Thus ${\rm exc}^{{\rm Abs}}(\sigma')={\rm exc}^{{\rm
Abs}}(\sigma)-2$. Also, we have ${\rm cyc}(\sigma') \equiv {\rm
cyc}(\sigma) \pmod 2$. Hence, the restriction of $a_n$ to $R_n^1$
is $q^2 b_{n-1}$.

In summary, we have:

$$a_n=(q-1)(a_{n-1}+qb_{n-1}),$$ and by symmetry :
$$b_n=(q-1)(b_{n-1}+qa_{n-1}).$$

Since $D_1=\{1\}$ we get $a_1=-1,b_1=-q$. Solving the above system of recursive equations yields:

$$P_{D_n}^{{\rm Abs}}(q,1,-1)=\frac{1}{2}(q-1)^{n-1}\sum_{\substack{k=0 \\ k {\textrm \quad {\rm even}}}}^n{{n \choose k} q^k}=
-\frac{1}{2}(q-1)^{n-1}((1+q)^n+(1-q)^n)
$$ as needed.

\section*{Acknowledgements}
We wish to thank Alex Lubotzky and Ron Livne. We also wish to
thank the Einstein Institute of Mathematics at the Hebrew
University for hosting their stays.


\begin{thebibliography}{LubW}
\bibitem{KZ}
G.\ Ksavrelof and J.\ Zeng, {\it Two involutions for signed excedance numbers},
Semi. Loth. Comb. {\bf 49} (2002/04), Art. B49e, 8 pp. (electronic).

\bibitem{Sta}
R.\ P.\  Stanley, {\it  Enumerative combinatorics}, Vol 1 and 2, Cambridge
University Press, 1997.

\end{thebibliography}
\end{document}